\documentclass[a4paper,12pt] {article} 
\usepackage{amsmath,amsfonts}
\usepackage{amsthm}

\DeclareMathSymbol{\varkappa}       {\mathord}{AMSb}{"7B}

\theoremstyle{plain}

\newtheorem*{thm1}{Theorem 1}
\newtheorem*{thm2}{Theorem 2}  
\newtheorem*{thm3}{Theorem 3}
\newtheorem*{thmA}{Theorem A}
\newtheorem*{thmB}{Theorem B}
\newtheorem*{thmC}{Theorem C}  
\newtheorem{lem}{Lemma}
\newtheorem*{pro1}{Proposition 1}
\newtheorem*{pro2}{Proposition 2}

\newtheorem*{cor2th3}{Corollary 2 of Theorem 3}
\newtheorem*{cor3th3}{Corollary 3 of Theorem 3}
\newtheorem*{cor1}{Corollary 1}

\theoremstyle{definition}
\newtheorem*{dfn}{Definition}

\newcommand{\nc}{\newcommand}
\newcommand{\rnc}{\renewcommand}

\rnc{\theequation}{\arabic{section}.\arabic{equation}}

\begin{document}

\title
{\textsc{\footnotesize To appear in Nagoya Mathematical Journal }\\
\vskip1cm
Vector semi-Fredholm \\
Toeplitz operators    \\
and mean winding numbers}

%
%

\maketitle

\footnotetext[1]{ This work was possible due to the support of the
Ministry of Science and Technology of Spain under the Ram\'on y
Cajal Programme (2002), the FEDER and the MEC grants MTM2004-03822
and MTM2005-08359-C03-01.}


\begin{abstract}
For a continuous nonvanishing complex-valued function $g$ on the real line,
several notions of a mean winding number are
introduced. We give necessary conditions for a Toeplitz operator
with matrix-valued symbol $G$ to be semi-Fredholm in terms of mean
winding numbers of $\det G$. The matrix function $G$ is assumed to
be continuous on the real line, and no other apriori assumptions
on it are made.
\end{abstract}




%


\nc\wh {\widehat}
\nc\ds\displaystyle
\rnc\Im {\operatorname{Im}}

\nc\beqn{\begin{equation}}
\nc\neqn{\end{equation}}
\nc{\beqnay}{\begin{eqnarray}}   \nc{\neqnay}{\end{eqnarray}}
\nc{\beqnays}{\begin{eqnarray*}} \nc{\neqnays}{\end{eqnarray*}}
\nc{\barr}{\begin{array}}        \nc{\narr}{\end{array}}
\nc\nn{\nonumber}

\nc{\lb}{\label}
\nc{\sn}{\section}               \nc{\ssn}{\subsection}
\nc{\nin}{\noindent}             \nc{\nl}{\newline}
\rnc{\theequation}{\arabic{section}.\arabic{equation}}
\rnc{\phi}{\varphi}
\rnc{\kappa}{\varkappa}
\rnc{\.}{.}
\nc{\df}[2]{{{\displaystyle#1}\over{\displaystyle#2}}}
\nc{\fr}{\frac}
\nc{\lra}{\longrightarrow}
\nc\rank{{\rm rank}\,}
\nc{\ind}{{\rm ind\,}}
\nc{\Ker}{{\rm Ker \,}}
\nc{\sign}{{\rm sign\,}}
\nc{\alg}{{\rm alg\,}}
\nc{\wind}{{\rm wind\,}}
\nc{\dist}{{\rm dist\,}}
\nc{\supp}{{\rm supp\,}}
\nc{\clos}{{\rm clos\,}}
\nc\ovl{\overline}
\rnc\Bbb{\mathbb}
\rnc\phi{\varphi}  
\nc\BR{\Bbb R}
\nc\BC{\Bbb C}
\nc{\BMO}{{\rm BMO}}
\nc{\BMOR}{{\rm BMO_{\BR}}}
\nc{\BMOp}{{\rm BMO_{\BR}^+}}
\nc{\BMOm}{{\rm BMO_{\BR}^-}}
\nc{\VMO}{{\rm VMO}}
\nc{\PC}{{\rm PC}}
\nc{\QC}{{\rm QC}}
\nc{\PQC}{{\rm PQC}}
\nc\Ltwo{L^2(\BR)}
\nc\Ltwor{L^2_r(\BR)}
\nc\Htwor{{H^2_r}}
\nc\Htwomr{{H^2_{-,r}}}
\nc{\Range}{{\rm Range\,}}
\nc{\Coker}{{\rm Coker\,}}
\nc\Ind{{\rm Ind\,}}
\nc{\St}{{\rm St\,}}
\nc{\Lip}{\rm{Lip}\,}
\nc\const{\text{const}}
\nc\ulim{\overline\lim}
\nc\llim{\underline\lim}
\nc\uw{{\overline w}}
\nc\lw{{\underline w}\,}
\nc\utw{{\widetilde w}}
\nc\ltw{{\underset {\sim} w}}

\nc\vka{\varkappa}
\nc\osmall{\overline{\overline{o\,}}}
\nc\defn {\overset {\text {\rm def} }{=}}
\nc{\Pol}{{\rm Pol\,}}
\nc\sm{\setminus}
\nc\gews{\geq}
\nc\lews{\leq}
\nc\al{\alpha}
\nc\be{\beta}
\nc\ze{\zeta}
\nc\ga{\gamma}
\nc\de{\delta}
\nc\om{\omega}
\nc\eps{\varepsilon}
\nc\Ga{\Gamma}
\nc\De{\Delta}
\nc\La{\Lambda}
\nc\CB{{\cal B}}
\nc\CO{{\cal O}}
\nc\CP{{\cal P}}
\nc\CA{{\cal A}}
\nc\CG{{\cal G}}
\nc\cD{{\cal D}}
\nc\cV{{\cal V}}
\nc\cW{{\cal W}}
\nc\CU{{\cal U}}
\nc\CH {{\cal H}}
\nc\CL{{\cal L}}
\nc\BD{{\Bbb D}}
\nc\BZ{{\Bbb Z}}
\nc\BN{{\Bbb N}}
\nc\BT{{\Bbb T}}
\nc\Crr{C_{r\times r}}
\nc\Cb{C^b}
\nc\Cbrr{C^b_{r\times r}}
\nc\GCbrr{{\cal G}C^b_{r\times r}}
\nc\CbR{{C^b(\BR)}}
\nc\Hnty{H^\infty(\BC_+)}
\nc\Hntys{H^\infty_*(\BC_+)}
\nc\GHnty{{\cal G}H^\infty(\BC_+)}
\nc\Hntyrr{H^\infty_{r\times r}(\BC_+)}
\nc\GLntyrr{{\cal G}L^\infty_{r\times r}}
\nc\GHntyrr{{\cal G}H^\infty_{r\times r}(\BC_+)}
\nc\HntyrrBD{H^\infty_{r\times r}(\BD)}
\nc\Lnty{L^\infty} \nc\Lntyrr{L^\infty_{r\times r}}
\nc\Carr{C_{a,r\times r}} \nc\GCarr{{\cal G}C_{a,r\times r}}
\nc\Aalrr{A^\al_{r\times r}}
\nc\Arrloc{A^{\al,\rm loc}_{r\times r}}
\nc\Arr{A^{\al}_{r\times r}}
\nc\Aloc{A^{\al,\rm loc}}
\nc\Alocrr{A^{\al,\rm loc}_{r\times r}}
\nc\Lal{{\rm Lip}^\alpha}
\nc\Lalrr{{\rm Lip}^\alpha_{r\times r}}
\nc\Lloc{{\rm Lip}^{\al,\rm loc}}
\nc\Llocrr{{\rm Lip}^{\al,\rm loc}_{r\times r}}
\nc\Om{\Omega}
\nc\clcD{\closu\cD}
\nc\la{\lambda}
\nc\HoneR{{H^1_{\BR}}}

\section{Introduction and main result}


Let $\BC_+=\{z\in\BC: \Im z>0\}$ be the upper half-plane in the
complex plane $\BC$.
We recall that the
classical Hardy space $H^p(\BC_+)$ consists of analytic functions
$f$ in $\BC_+$ such that $\|f\| \defn  \Big(
\sup_{y>0}\int_{\BR}|f(x+iy)|^p\,dx \Big)^{1/p} $ is finite.
It is a Banach space for any $p$ as above.
The space $H^\infty(\BC_+)$ is defined as the Banach space
of bounded analytic functions in $\BC_+$.
We refer to the book \cite{Dur} for an account of the
theory of $H^p$ spaces of the upper half-plane and of the unit disc.
Functions in $H^p(\BC_+)$ have non-tangential boundary limit
values on $\BR$, which permits us to identify $H^p(\BC_+)$ with a
closed subspace of $L^p(\BR)$.
We put $H^p=H^p(\BC_+)$, $1\le p\le \infty$.

For any function space $\Psi$, we denote by
$\Psi_\BR$ the set of its real elements and by
$\Psi_r$, $\Psi_{r\times r}$, respectively,
the spaces of $r\times 1$ vector-valued functions and of
$r\times r$ matrix-valued functions with entries in $\Psi$.
If $\CA$ is a scalar or matrix functional algebra, we denote by
$\CG\CA$ the set of all its invertible elements.

Let the natural number $r$ be fixed and let $G\in L^\infty_{r\times r}(\BR)$.
The vector Toeplitz operator $T_G$ with the symbol $G$
acts on the vector Hardy space $\Htwor$ by the formula
\beqn
\label{TG}
T_Gx=P_+\big(G\cdot x\big), \quad x\in\Htwor,
\neqn
here $P_+$ is the orthogonal projection of
$\Ltwor$ onto its closed subspace $\Htwor$.

A (bounded linear) operator $K$ on a Banach space $B$ is called
normally solvable \cite{GK}, \cite{MikhPr} if its image is closed.
$K$ is called a $\Phi_+$ operator (a $\Phi_-$ operator) if it is
normally solvable and $\dim\Ker K<\infty$ ( $\dim\Coker K=\dim B /
\Range K<\infty$, respectively). We denote by $\Phi_\pm(B)$ these
classes of operators on $B$. Operators in $\Phi_+$ and $\Phi_-$
are called semi-Fredholm. Operators in $\Phi(B)=\Phi_-(B)\cap
\Phi_+(B)$ are called Fredholm.

The index of a semi-Fredholm operator is defined by
$$
\Ind K=\dim\Ker K -\dim\Coker K;
$$
its values are integers or $\pm\infty$.
A semi-Fredholm operator is Fredholm if and only if its index is finite.

Fredholm and semi-Fredholm operators have several important
properties. For instance, the product of two $\Phi_\pm$ operators
is again a $\Phi_\pm$ operator, and the formula
$\Ind(K_1K_2)=\Ind(K_1)+\Ind(K_2)$ holds for $K_1$, $K_2$ in
$\Phi_+$ or in $\Phi_-$. We refer to \cite{GK}, \cite{MikhPr} for
detailed expositions of the theory of these classes and for
applications.


We put $\Cb=\Cb(\BR)$ to be the Banach space of all continuous
bounded functions on $\BR$ with the supremum norm. Our
paper is devoted to finding necessary conditions for
semi-Fredholmness and Fredholmness of $T_G$ for the case when $G$
is an $r\times r$ matrix function whose entries are in $\Cb$. Such
questions appear naturally in connection with the Riemann--Hilbert
problem on the real line. This problem appears in many different
situations, such as various problems in mechanics of continuous
media and hydrodynamics \cite{Begehr}, \cite{MonSem},
\cite{AntPopYats}, \cite{Beg_Dai}, \cite{ShargTol}, inverse
scattering method for integrable equations \cite{AblCla}, linear
control theory of systems with delays \cite{CallWi}, convolution
equations and systems on finite intervals (see \cite{BKSpi},
\cite{KSpi_Semifr}) and others. The case of infinite index often
appears in these applications.

First we quote the following well-known result.

\begin{thmA} [see  \cite{BSlb}]
The condition $\det G\ge\eps>0$ is necessary for $T_G$ to be
semi-Fredholm.
\end{thmA}

We will always assume this condition to be fulfilled.

For a function $G\in \Cbrr$ which has limits
at $\pm\infty$,
$T_G$ is semi-Fredholm iff it is Fredholm, and
a complete criterion for it is known
(see \cite{BSlb} or \cite{BKSpi}).
In a particular case, when
$G(-\infty)=G(+\infty)$,
$T_G$ is Fredholm if and only if $|\det G|\ge\eps>0$ on
$\BR$, and
\beqn
\label{ind-wind}
\ind T_G=-\wind\det G,
\neqn
where $\wind $ stands for  the winding number
(around the origin).
So our main concern is about symbols that
have no limits either at $-\infty$ or at $+\infty$.

Let $\BMO=\BMO_\BR$ be the space of real-valued functions on $\BR$
of bounded mean oscillation. We recall that $\BMO$
consists of those locally integrable functions $f$ on $\BR$ that satisfy
\beqn
\label{defBMO}
\sup_J\frac 1 {|J|}\int \big|f-f_J  \big|\le C,
\neqn
where the supremum is taken over all finite subintervals $J$ of
the real line and $f_J=\frac 1 {|J|}\int_J f$ is the mean of $f$
on the interval $J$. It is known \cite{GR} that
if there exist
a constant $C$ and
arbitrary real numbers $f_J$ such that
(\ref{defBMO}) holds for any finite interval
$J$ in $\BR$, then $f$ belongs to $\BMOR$.
We refer to \cite{GR} for an exposition of the theory of these
spaces.

Let $C_+(\BR)$ be the class of real continuous
(nonstrictly) increasing functions on $\BR$, and put
\beqnay
\BMOp= \big\{ u+v: \quad u\in\BMOR, \,v\in C_+(\BR) \big\},
\nonumber                  \\
\BMOm=
\big\{
u-v: \quad u\in\BMOR, \, v\in C_+(\BR)
\big\}. \nonumber
\neqnay
The main result of \S2 is as follows.


\begin{thm1}
{\it Suppose that $G\in \Cbrr$.

(1) If $T_G\in\Phi_\pm(\Htwor)$, then
$\arg \det
G\in\BMO^\pm_\BR$.


(2) If $T_G\in\Phi(\Htwor)$, then $\arg \det G\in\BMOR$. }
\end{thm1}

In \S3, we introduce a system of mean winding numbers of
 $\det G$ and formulate and prove
Theorems 2 and 3 (they will follow from Theorem 1). In
\S4, we discuss some unresolved questions, related with our results.

Our principal motivation comes from the control theory. In a
problem about the complete controllability of delay equations
it turned out to be necessary to estimate the number
\beqn
\inf\, \Big\{
\tau\in \BR: \quad T_{e^{-i\tau x}G(x)} \text{ is onto}
\Big\}\,\defn\,\beta(G)
\label{betaG}
\neqn
in terms of some computable characteristics of a matrix function $G\in \GCbrr(\BR)$.
The number $\beta(G)$ has a meaning
of the least time of complete controllability.
Theorems 2 and 3 permitted us to give
a good estimate of this number.
These results were obtained
jointly by the author and Sjoerd Lunel
and will be published elsewhere.

A great part of the recent book \cite{DybGru} by Dybin and Grudsky
treats scalar and matrix functions that are continuous on the real
line. This book summarized
(and generalized) earlier work by these authors.
Several novel
tools are used, such as the notion of a $u$-periodic function,
where $u$ is an inner function on $\BC_+$, continuous on the real
line. Another tools are a construction of an inner function whose
argument models an arbitrarily given increasing continuous
function and the notion of a generalized factorization with
infinite index. These hard analysis
tools permitted the authors to give a sufficient condition for
semi-fredholmness (see \cite{DybGru}, Theorem 5.10. By applying
this result, Dybin and Grudsky
get complete answers in cases of whirls at
$\pm \infty$ with different asymptotic, such as power, logarithmic
or exponential.

Earlier work on whirled symbols include the works by Govorov
\cite{Gov}, Ostrovsky \cite{Ostrovs}, Monakhov, Semenko (see the
book \cite{MonSem}) and others; the approach of these authors is
based on the theory of analytic functions of completely regular
growth. In various works, the behavior of the property of
Fredholmness under an orientation preserving homeomorphism of
$\BR$ have been studied, see \cite{BGRa}, \cite{DybGru},
\cite{BGSpi2} and others.

Various mean winding numbers were introduced in the work
by Sarason \cite{Sar} for symbols in $\QC$ and by Power \cite{Pow}
for slowly oscillating symbols.
For symbols of these classes, these
mean winding numbers allow one to
formulate nice complete criteria for a Toeplitz operator
to be Fredholm or semi-Fredholm. We remark that a wider $C^*$-algebra of
slowly oscillating functions was considered in a recent paper
by Sarason \cite{Sar3}, where the maximal ideal space
of this algebra was studied.

Necessary and sufficient conditions for a
Toeplitz operator to be Fredholm and semi-Fredholm are also known
if $G$ belongs to various algebras of symbols. For instance,
classes $\PC$ of piecewise continuous symbols, $\QC=L^\infty\cap
\VMO$ of quasicontinuous symbols, and $\PQC=\alg(\PC,\QC)$ have
been studied both in scalar and matrix case.

Another well-studied cases are that of almost periodic and
semi-almost periodic symbols. For matrix symbols of these types, a
great breakthrough has been done recently by B\"ottcher, Karlovich
and Spitkovsky, see \cite{BKSpi}. Among other things,
generalizations of the index formula (\ref{ind-wind}) are known
for these cases (see \cite{BSlb}, \cite{Nik}, \cite{NikOFS}).
We refer to \cite{BKSpi2003} for an alternative approach.
In  \cite{Abr, BasFernKarl, BGRa},
other classes of symbols are studied. In \cite{BKS04},
a Fredholm criterion and an index formula are given for vector Toeplitz
operators, whose (matrix) symbols belong to the Banach algebra, generated
by semi-almost periodic matrix functions and slowly
oscillating matrix functions.
See \cite{LSpi} for a connection
with the factorization and the Riemann--Hilbert problem.

%
%
%

For symbols in $\Cbrr(\BR)$ with no other assumptions, our
knowledge is much less complete. We refer to Subsections 2.26 and
4.73 in \cite{BSlb} and to \cite{BG} for several relevant results.
The criterion for surjectivity of a Toeplitz operator
with a nontrivial kernel, given in \cite{HrSarSeip}, can also be
reformulated as a criterion for a Toeplitz operator to belong to
$\Phi_+\sm\Phi$.
Some additional comments will be given at the end of the article.

We refer to \cite{LSpi}, \cite{Nik}, \cite{NikOFS}, \cite{BKSpi} for
systematic expositions of the spectral theory of Toeplitz operators.


It is worth to note that recently, Toeplitz operators with symbols
like ours have been appeared in papers by Baranov, Havin, Makarov,
Mashreghi, Poltoratsky and others in relation with the
Beurling--Malliavin theorem, bases in de Branges spaces and
related topics (see \cite{HavMashr1-2}, \cite{MakPolt2},
\cite{BarHav} and references therein). It seems that the ideas and
methods of these papers can be applied to achieve a better
understanding of semi-Fredholm Toeplitz operators with continuous
symbols at least in the case of scalar symbol $G$.



\textbf{Acknowledgements.} The author expresses his gratitude to
M. Gamal and I. Spitkovsky for valuable comments.


\section{Proof of Theorem 1}
\setcounter{equation}{0}

First we need some facts and definitions.

Let $0<\al<1$. We put
$$
\Lloc=\{f\in \Cb(\BR):f|J\in \Lal(J)\quad \forall J\}; \nonumber
$$
here $J$ runs over all compact intervals in $\BR$ and
$\Lal(J)$ is the H\"older--Lipschitz class on $J$ with the exponent $\al$.
Next, we will need the classes
\begin{align}
C_a&=\{f\in C(\clos \BC_+):f|\BC_+\in H^\infty\}, \nonumber \\
\Aloc&=\{f\in C_a:f|\BR\in \Lloc\}. \nonumber
\end{align}
A function $f$ in $\Llocrr$ or in $\Alocrr$ is invertible
if and only if $|\det f|>\eps>0$ on $\BR$ (or on
$\clos\BC_+$, respectively).
%
%
Recall that a function $g$ in $\Hnty$ is called {\it inner} if its
modulus is equal to one a.e. on $\BR\,$. The function $g$ is
called {\it outer} if it has a form $g(z)=\exp(u(z)+iv(z))$,
$$
(u+iv)(z)=\frac1{\pi i}\int_{-\infty}^{\infty} \bigg[
\frac1{t-z}-\frac t{1+t^2} \bigg] \log
k(t)\,dt+is,
$$
where $k>0$ a.e. on $\BR$, $\log k\in L^1(\BR)$, and $s$ is a
real constant. We assume $u$ and $v$ to be real-valued and harmonic in $\BC_+$.
These functions have boundary limit values a.e. on
$\BR$, which satisfy $u|\BR=\log k$ a.e. and $v|\BR=
\CH(u|\BR)$, where $\cal H$ is the Hilbert transform on $\BR$.

Each function $g$ in $\GHnty$ is outer; in this case
$\log k\in L^\infty(\BR)$. We refer to
\cite{Dur}, \cite{GR} for all these (classical) facts.

For each function $g$ in $\GHnty$, $\arg g(z)=s+v(z)$ is
well-defined on $\BC_+$ (up to an additive constant
$2\pi n$). We also see that the function $\arg g(z)$ has boundary
limit values a.e. on $\BR$, which will be denoted as $\arg g(x)$,
$x\in\BR$.

\nin{\bf Definition.} \enspace We define the class
$\Hntys$ as the set of functions
$f\in\Hnty$ that have the form
$$f=g\cdot h,
$$
where $g\in\GHnty$ and $h$ is inner in $\BC_+$ and has a
continuous extension to $\BR$.

A function $h$ is inner of the above type if and only if
it has the form
\beqn
\label{inner cont on R}
h(z)=Ce^{iaz}
\prod_j
\frac
{|z_j^2+1|} {z_j^2+1} \,
 \frac {z-z_j} {z-\bar z_j}, \qquad
z\in\BC_+,
\neqn
where $|C|=1$, $a>0$, and $z_j\in\BC_+$, $|z_j|\to\infty$. Take
any positive continuous function $y=\psi(x)$ on $\BR$ such that
the subgraph $\Ga_\psi= \{(x+iy):\, 0<y<\psi(x)\}\subset\BC_+$
does not contain the zeros $z_j$ of $h$. Then $\arg h(z)$ is
well-defined and continuous on $\Ga_\psi\cup\BR$.

\nin{\bf Definition.} \enspace
Let $f\in\Hntys$, and let $g$, $h$, $\Ga_\psi$ be as above.
We define the argument $\arg f$ on $\Ga_\psi\cup \BR$ by
$$
\arg f=\arg g+\arg h.
$$

\nin So for $f\in\Hntys$, the argument $\arg f$ is well-defined on
$\Ga_\psi$ (up to adding $2\pi n$, $n\in \BZ$). It is continuous
on $\Ga_\psi$ and its values on $\BR$
exist almost everywhere in the sense of nontangential limits.

\begin{pro1}
For any $f\in\Hntys$, $\arg f\in\BMOp$.
\end{pro1}
\begin{proof}
For any $f=g\cdot h\in\Hntys$ as above, $\arg g\in\BMOR$ and $\arg
h$ is a continuous increasing function.
\end{proof}

\begin{lem}
\label{lem Hntys}
Let $f\in H^\infty$. Then $f\in H^\infty_*$ if and only if there is a
positive function $\psi\in C(\BR)$ and some $\eps>0$
such that $|f|>\eps$ on the subgraph $\Ga_\psi$.
\end{lem}

\begin{proof}
If $f\in H^\infty_*$, then it is clear that $f$ satisfies the above property.
Conversely, suppose $|f|>\eps>0$ on $\Ga_\psi$, for a certain
positive function $\psi\in C(\BR)$. Let $f=h\cdot g$ be the inner - outer
factorization of $f$, then $g\in {\cal G}H^\infty$. It follows that
the inner function $h=f/g$ satisfies an inequality $|h|>\eps_1>0$
on $\Ga_\psi$, and consequently, it has a form (\ref{inner cont on R}), see
\cite[Chapter 3]{Nik}.
\end{proof}

In many works on Toeplitz operators, the
unit disc $\BD=\{z\in\BC: |z|<1\}$ instead of
the upper half-plane $\BC_+$ is considered. If $G\in
L^\infty_{r\times r}(\BT)$, where $\BT=\partial \BD$ is the unit
circle, then the same formula
(\ref{TG}) defines a Toeplitz operator $\wh T_G$ on $H^2_r(\BD)$
(in this setting, $P_+$ stands for the orthogonal projection of $L^2_r(\BT)$
onto the vector Hardy space $H^2_r(\BD)$).
Let
\beqn
\varphi(z)=\frac{z-i}{z+i}
\neqn
be the conformal mapping of $\BC_+$ onto the unit disc $\BD$.
The formula
\beqn
T_G=W \wh T_{G\circ \phi} W^{-1},
\neqn
where $W: \Htwor(\BD)\to\Htwor$ is the unitary isomorphism, given by
$$
(Wf)(z)=\pi^{-1/2}\,(z+i)^{-1}
\big(f\circ\varphi\big)(z)
$$
shows that each vector Toeplitz operator on $\BC_+$ is unitarily
equivalent to a vector Toeplitz operator on $\BD$, and vice versa, so there
is no difference in the study of Toeplitz operators in these two settings.
The symbols on $\BT$ that correspond to symbols in
$\Cbrr$ by means of this construction
have the only discontinuity at the point $1$.


Notice first of all that each function $G$ in $\GLntyrr(\BT)$
factors (in an essentially unique way) as $G=UG_e$, where
$G_e\in\CG\HntyrrBD$ and $U$ is unitary-valued on $\BT$. Then
$T_G=T_UT_{G_e}$, and $T_{G_e}$ is invertible, so that
Fredholmness or semi-Fredholmness of $T_G$ is equivalent to the
corresponding property of $T_U$. For unitary symbols, the
following results hold, see \cite{ClaGo}, \cite{BSlb}.

\begin{thmB} 
Let $U\in\GLntyrr(\BR)$ be unitary-valued. Then

(i) $T_U$ is left-invertible
if and only if $\dist(U, \,\Hntyrr)<1$.

(ii) $T_U$ is invertible if and only if $\dist(U,\, \GHntyrr)<1$.
\end{thmB}

\begin{thmC}   
Let $U\in\GLntyrr(\BT)$ be unitary-valued. Then

(i) $\wh T_U\in \Phi_+$ if and only if $\dist(U,\, \Crr(\BT)+\HntyrrBD)<1$.

(ii) $\wh T_U\in \Phi$ if and only if $\dist(U,\, \CG\big(\Crr(\BT)+\HntyrrBD )\big)<1$.

\end{thmC}
We refer to \cite{BSlb}, Section 4.38 for the connection with Fredholmness.

We will also make use of the following properties.

\begin{pro2}
\label{pro2}
\begin{itemize}
\item[(1)] Each selfadjoint matrix function
$K\in\Lntyrr(\BR)$ such that
$K(x)\ge\eps I>0$ on $\BR$ has a factorization
$K(x)=G^*_e(x)G_e(x)$ on $\BR$, where
$G_e\in\GHntyrr$.
This factorization is unique up to multiplying
$G_e$ on the left by a constant unitary matrix.

\item[(2)] If the matrix $K$ (as above) satisfies additionally
$K\circ \phi\in\Lalrr(\BT)$, then  $G_e\circ\phi\in {\cal G}\Arr(\clos \BD)$; here
$$
A^\al(\clos \BD)=
\{
f\in C(\clos \BD): f|\BD\in H^\infty(\BD), \;
f|\BT\in{\rm Lip\,}^\al(\BT)
\}.
$$
\end{itemize}
\end{pro2}

For the property (1), see \cite{LSpi}, Theorem 7.9 and
\cite{SzNF}, Proposition 7.1. The proof of
(2) is contained in \cite{Shmu}, \cite{BSlb}.

\begin{lem} 
\label{lm185}
Let $G\in\Lntyrr(\BR)$. Then $T_G\in\Phi_+$ if and only if $T_{\phi^nG}$
is left invertible for some integer $n\gews0$.
\end{lem}

\begin{proof} It is more transparent
to work with $\Htwor(\BD)$ instead of
$\Htwor$. Suppose $G=G(z)\in\Lntyrr(\BT)$ and $\wh T_G\in\Phi_+$; we have to
check that there is some
integer $n\gews0$ such that $\wh T_{z^nG}$ is left invertible. By
the assumption, the kernel
$\Ker\wh T_G$ is finite dimensional; let $x_1, \dots x_m\in \Htwor(\BD)$
be its basis. Put
$$
L_n=\{(c_1,\dots c_m)\in\BC^m: G\cdot \sum_j c_jx_j\in z^{-n}\Htwomr\},
\quad n\gews 0,
$$
where $H^2_{-,r}=L^2_r(\BT)\ominus \Htwor(\BD)$.
Then
$
\BC^m=L_0\supset L_1\supset \dots \supset L_n\supset \dots\qquad.
$
Since $\bigcap_0^\infty L_k=0$, one has $L_n=0$ for some $n\gews0$.
If $x\in \Ker \wh T_{z^nG(z)}$, then $x=\sum_{j=1}^m c_jx_j$ for some
coefficients $c_j$ and $z^nGx\in \Htwomr$, which implies that
$c_1=\dots=c_m=0$.
Hence $\Ker \wh T_{z^nG(z)}=0$.
Since $\wh T_{z^nG}=\wh T_G \wh T_{z^n}$ is a $\Phi_+$ operator
with trivial kernel, it follows that it is
left invertible.
\end{proof}

\begin{lem}
\label{intersBMO}
1) Suppose that $u_1, u_2$
are real increasing functions on $\BR$ and $u:=u_1+u_2\in \BMOR$. Then
$u_1,u_2\in \BMOR$.

2) $\BMOm\cap \BMOp=\BMOR$.
\end{lem}

\begin{proof} 1) For any finite interval $J\subset \BR$, one can find a point
$c=c_J\in J$ such that $u(x)\lews u_J$ for $x<c_J$ and
$u(x)\gews u_J$ for $x>c_J$. There exist numbers
$\al_{1J}$, $\al_{2J}$ such that
$u_k(c_J-0)\lews \al_{kJ}\lews u_k(c_J+0)$ for $k=1,2$ and
$\al_{1J}+\al_{2J}=u_J$. Then for any subinterval
$J\subset \BR$,
$$
\int_J |u_1(x)-\al_{1J}|\, dx
+\int_J |u_2(x)-\al_{2J}|\, dx
=\int_J |u(x)-u_J|\, dx \le C|J|,
$$
where $C=\|u\|_{\BMOR}$.
It follows that $u_1,u_2\in\BMOR$.

2) If $h=w_1-v_1=w_2+v_2\in\BMOm\cap\BMOp$, where
$w_1,w_2\in \BMOR$ and
$v_1,v_2\in C_+(\BR)$, then by part 1),
$v_1,v_2\in \BMOR$ because
$v_1+v_2\in \BMOR$.
\end{proof}

The next lemma is not new; in fact, Spitkovsky gives in
\cite[Theorem 2]{Spi2} a more general result. We will give a proof
for completeness.

\begin{lem}
\label{lm186}
\enspace    
Suppose that $J$ is a finite open interval
on the real line, $F, G\in\GHntyrr$, and
$F^*F=G^*G$ a.e. on $J$.
Then there exists a neighbourhood $\cal W$ of
$J$ in $\BC$ and a  bounded  analytic $r\times r$ matrix function
$V$ on $\cal W$
such that $F=VG$ on $\cal W$ (and a.e. on $J$)
and $V$ is unitary valued on $J$.
\end{lem}

\begin{proof}
Put $V=FG^{-1}$, then $F=VG$ on $\BC_+$ and a.e. on $\BR$ and
$V$ is unitary on $J$.  We apply the symmetry
principle to $V$. Since $V\in\GHntyrr$, it is easy to prove that
$\widetilde V(z)=V^{*-1}(\bar z)$ is an analytic continuation of $V$ onto
the lower half-plane through the arc $J$.
\end{proof}

\begin{lem}
\label{lm196}
Every matrix function $G\in\Llocrr$
such that
$\inf_\BR |\det G|>0$
has a factorization $G=UG_e$, where $G_e$, $G_e^{-1}\in \Alocrr$
and $U\in\Llocrr$ is unitary valued.
\end{lem}

\begin{proof}
Put $K(x)=G^*(x)G(x)$, then
$K(x)\ge\eps_1 I > 0$ on $\BR$. By the above property (1),
$K$ can be factorized as $K(x)=G_e^*(x)G_e(x)$,
where $G_e\in\GHntyrr$.
Hence $G=UG_e$, where $U\in\Lntyrr$.

Consider a sequence of
matrix functions $K_n$ such that
$K_n(x)=K(x)$ on $[-n,n]$, $K_n(x)\ge \eps_1 I >0$ on $\BR$ and
$K_n\circ \phi$ are Lipschitz on $\BT$. By property (2),
we arrive at functions
$G_{ne}\in\GHntyrr$ such that
$G_{ne}\circ \phi\in\Arr(\BD)$ and
$K_n=G_{ne}^*G_{ne}$ on $\BR$. By
Lemma \ref{lm186}, $G_e=V_nG_{ne}$ on $(-n,n)$,
where $V_n$ are unitary on $(-n,n)$ and
analytic in neighbourhoods of these
intervals. It follows that $G_e\in \Alocrr$.
Therefore $U\in\Llocrr$.
\end{proof}


\begin{lem} 
\label{lm188}
Suppose that $H\in\Cbrr(\BR)$ and $\Psi\in \Hnty$.
Then for any finite interval $L$ on the
real line we have
\begin{equation}
\limsup_{y\to0+}\|\Psi(\cdot+iy)-H(\cdot)\|_{\Lntyrr(L)}\lews
\|\Psi-H\|_\infty;                       \label{AA}
\end{equation}
here
$\|\Psi-H\|_\infty=\|\Psi-H\|_{\Lntyrr(\,\BR\,)}$.
\end{lem}

\begin{proof}
Denote by $H(z)$, $z\in \BC_+$, the harmonic extension of $H$
by means of the Poisson formula. Then for any $y>0$,
$$
\|\Psi(\cdot +iy)-H(\cdot+iy)\|_{\Lntyrr(L)}\lews \|\Phi-H\|_\infty.
$$
Since $H(x+iy)\to H(x)$ as $y\to0+$ uniformly on
compact subsets of the real line, the result follows.
\end{proof}


\begin{proof}[Proof of Theorem 1]
We prove part (1).
Suppose
$G\in \Cbrr$ and
$T_G\in\Phi_+$. We have to prove that
$\arg\det G\in\BMOp$.
By Lemma \ref{lm185}, there is some
$k>0$ such that $T_{G_1}$ is left invertible, where
$G_1=\phi^kG$. Since $\arg\det G=
\arg\det G_1-kr\arg\phi$ and
$\arg\phi\in L^\infty(\BR)\subset\BMOR$,
we have only to prove that
$\arg\det G_1$ is in $\BMOp$.
Let $\|T_{G_1}x\|\gews\eps\|x\|$, $x\in H^2_r$,
where $\eps\gews0$, then
for any $G_2$ with $\|G_1-G_2\|_\infty<\eps$,
$T_{G_2}$ is also left invertible. Take
$G_2=G_1+R$ such that
$G_2\in\Llocrr$ and
$R\in C^b_{r\times r}$ has a  small norm
$\|R\|_\infty$:
$\|R\|_\infty<\eps'<\eps$, where $\eps'$ has to be chosen.
Since
$$
\arg\det G_2=\arg\det G_1
+\arg\det (I+G_1^{-1}R),
$$
it follows that
$\arg\det G_2-\arg\det G_1\in L^\infty(\BR)$
if we assume that
$\eps'\cdot\|G_1^{-1}\|_\infty<1$.
So it suffices to consider $G_2$ instead of $G$.

By Lemma \ref{lm196}, we have a factorization
$G_2=UG_{2e}$, where
$U\in\Llocrr$ is unitary valued and
$G_{2e}\in \CG\Alocrr$.
Then
\beqn
T_{G_2}= T_U\, T_{G_{2e}}.
\nn
\neqn
Since $T_{G_{2e}}^{-1}=T_{G_{2e}^{-1}}$,
we conclude that $T_U$ is
left invertible. We apply Theorem B
and arrive at a function $F\in\Hntyrr$
with $\|U-F\|_\infty<1-\eps_0<1$. Put
$F_y(x)=F(x+iy)$, $y>0$,
$L=L_\rho=[-\rho, \rho]$, where
$\rho>0$. By Lemma \ref{lm188},
\begin{equation}
\|I-U(x)^{-1}F_y(x)\|_{\Lntyrr(L_\rho)}=
\big\|U(x)-F_y(x)\big\|_{\Lntyrr(L_\rho)}< 1-\eps_0      \label{AD}
\end{equation}
for $x\in L_\rho$, $y\in(0,\de)$, where
$\de=\de(\rho)>0$.
It follows, in particular, that
there is a graph $y=\psi(x)$ of a
positive function $\psi\in C^b(\BR)$
such that
$$
\|I-U(x)^{-1}F(x+iy)\|<1-\eps_0 \quad \text{for }x+iy\in\Ga_\psi.
$$
It follows that $\arg\det F$ is well defined on $\Ga_\psi$.
By Lemma \ref{lem Hntys}, $\det F$ belongs to $\Hntys$.

One can define a continuous branch of
$\arg\det\big(U(x)^{-1}F(x+iy)\big)$
for
\linebreak     
$x+iy\in\Ga_\psi$ so that
$|\arg\det\big(U(x)^{-1}F(x+iy)\big)|<r \pi/ 2$.
Therefore there is a continuous branch of
$\arg\det  F(x+iy)$, $x+iy\in\Ga_\psi$ such that its limit values satisfy
$$
\big|
\arg\det F(x)-\arg\det U(x)
\big|
\lews \frac {r \pi} 2 \qquad\text{ a.e. on }\BR.
$$
By Proposition 1,  $\arg\det F\in\BMOp$. Hence $\arg\det U\in\BMOp$.
Since $G_{2e}\in\CG\Arrloc(\BC_+)$,
it follows that
$\det G_{2e}\in \GHnty$,
so that
$\arg\det G_{2e}\in\BMOR$.
Finally, we deduce from the formula
$$
\arg\det G_2=
\arg\det U+\arg\det G_{2e}
$$
that $\arg\det G_2\in \BMOp$.

The case when $T_G\in\Phi_-$ is obtained by considering
$G^*$ instead of $G$.
The assertion (2) follows from (1) and Lemma \ref{intersBMO}.
\end{proof}

I. M. Spitkovsky communicated to the author an outline
of an alternative proof of Theorem 1, which is based on some properties of the
transplantation of the algebra $H^\infty(\BD)+C(\BT)$ to the real line.

\section{Mean winding numbers}
\setcounter{equation}{0}

Let $H^1_{\BR}$ be the real Hardy space,
\beqn
H^1_{\BR}=\big\{
u\in L^1_\BR(\BR): \CH u\in L^1_\BR(\BR)
\big\}.
\neqn
We put $\|u\|_{H^1_{\BR}}=\|u\|_{L^1}+\|\CH u\|_{L^1}$.

Consider the cone
$$
\Pi=\{ \eta\in\HoneR: \quad \eta  \text{
has a compact support on }\BR,  \quad
\int_{-\infty}^x\eta\le0\quad \forall x\in\BR\}.
$$

\begin{thm2}
Let $G$ be a $r\times r$
matrix function in $\Cbrr$.

(1) If $T_G\in\Phi_-(\Htwor)$,
there is a constant $C>0$ such that
for any $\eta$ in ~$\Pi$,
$$
\int_\BR\eta(x)
\big(\arg\det G\big)(x)\,dx\lews C\|\eta\|_{\HoneR}.
$$

(2) If $T_G\in\Phi_+(\Htwor)$
there is a constant $C>0$ such that for any $\eta$ in ~$\Pi$,
$$
\int_\BR\eta(x)
\big(\arg\det G\big)(x)\,dx\gews -C\|\eta\|_{\HoneR}.
$$
\end{thm2}

It is well known that $\int_\BR\eta=0$ for any function $\eta$ in
$\HoneR$, see \cite{GR}, Chapter III.
Hence
the above integrals do not depend on the
additive constant in $\arg\det G$.

As a consequence, we obtain that if $T_G\in\Phi(\Htwor)$, then
\beqn
\bigg|\int_\BR\eta(x)
\big(\arg\det G\big)(x)\,dx\bigg|\lews C\|\eta\|_{\HoneR},
\qquad \eta\in \Pi.
\nn
\neqn
In the scalar case, this inequality  follows from
the Widom-Devinatz theorem (Theorem B),
together with the Fefferman duality theorem,
and takes place for all
$\eta \in \HoneR$ (the integral is to be
understood in the sense of the
the duality $\HoneR$ -- $\BMOR$).

\begin{dfn}
Let $\eta\in\Pi$, $\eta\not\equiv0$ be fixed, and let $G\in\Cbrr$.
Define the upper and the lower \textit{mean winding numbers} of
$\det G$ (associated with $\eta$) by
\begin{align*}
\displaystyle
\uw_\eta(G)&=
\underset{T\to+\infty} {\ulim} \;
\sup_{y\in\BR}
\frac1T\int_{\BR}\eta\big(\frac{x-y}T\big)\cdot\arg\det G(x)\,dx, \\
\lw_\eta(G)&=
\underset{T\to+\infty} {\llim} \;
\inf_{y\in\BR}\,
\frac1T\int_{\BR}\eta\big(\frac{x-y}T\big)\cdot\arg\det G(x)\,dx.
\end{align*}
\end{dfn}


\begin{thm3}
(1) If $T_G\in\Phi_+(\Htwor)$, then $\lw{}_\eta(G)\ne-\infty$;

(2) If $T_G\in\Phi_-(\Htwor)$, then $\uw_\eta(G) \ne+\infty$.
\end{thm3}

One can also define simpler characteristics
\begin{equation*}
\displaystyle
\utw_\eta(G)=
\underset{T\to+\infty} {\ulim} \;
\frac1T\int_{\BR}\eta\big(\frac xT\big)\cdot\arg\det G(x)\,dx
\end{equation*}
and the number
$\ltw{}_\eta(G)$, defined as the corresponding lower limit.
One has
$\lw_\eta(G)\le
\ltw{}_\eta(G)\le
\utw_\eta(G)\le
\uw_\eta(G)$, so that Theorem 3 implies the same
assertions for
$\ltw{}_\eta(G)$, $\utw_\eta(G)$.

Consider a scalar $G\in\CbR$, $|G|>\eps>0$ on $\BR$. If $\arg
G$ has finite limits on $\pm\infty$, then
$\utw_\eta(G)=\ltw{}_\eta(G)= K\cdot \arg G\big|_{-\infty}^{+\infty}$,
where $K=\int_0^{+\infty}\eta(x)\,dx$.
One also has
$\uw_\eta(G)=L\cdot \big(\arg G\big|_{-\infty}^{+\infty}\big)_+$,
$\lw_\eta(G)=L\cdot \big(\arg G\big|_{-\infty}^{+\infty}\big)_-$,
where $L=\sup_{y\in\BR}\int_y^{+\infty}\eta$,
$y_+=\max(y,0)$, $y_-=\min(y,0)$.
So in this case all these
winding numbers have a simple sense.
For these symbols, each of the
conditions $T_G\in\Phi_-(\Htwor)$, $T_G\in\Phi_+(\Htwor)$,
$T_G\in\Phi(\Htwor)$ is equivalent to the requirement
$(\arg G)\big|_{-\infty}^{+\infty} \ne\pm\pi, \pm 3\pi, \pm 5\pi$, etc.
(see, for instance, \cite{BSlb} or \cite{GK}, Ch. 9).

\begin{cor1}
Let $\al>0$, and define generalized winding numbers
\begin{align*}
\displaystyle
&\uw_{\eta,\al}(G)=
\underset
{T\to+\infty}
{\ulim}
\sup_{y\in\BR}\;
\frac1{T^{1+\al}}
\int_{\BR}\eta\big(\frac{x-y}T\big)\cdot(\arg\det G)(x)\,dx, \\
&\lw_{\eta,\al}(G)=
\underset
{T\to+\infty}
{\llim}
\inf_{y\in\BR}\;
\frac1{T^{1+\al}}
\int_{\BR} \eta\big(\frac{x-y}T\big) \cdot(\arg\det G)(x)\,dx.
\end{align*}

(1) If $T_G\in\Phi_+(\Htwor)$, then $\lw_{\eta,\al}(G)\gews0$;

(2) If $T_G\in\Phi_-(\Htwor)$, then $\uw_{\eta,\al}(G)\lews0$.
\end{cor1}

This follows immediately from Theorem 3.                      \hfill $\Box$

In particular, the function $\eta_\al=\frac{1+\al}2
\big(\chi_{[0,1]}-\chi_{[-1,0]}\big)$ is in $\Pi$. The
corresponding upper winding number is given by
\beqn
\uw_\al(G)=
\underset
{T\to+\infty}
{\ulim}\;
\frac{1+\al}{2T^{1+\al}}\;
\sup_{y\in\BR}\;
\bigg[
\int_y^{T+y}
-\int_{y-T}^y \bigg] \;  \arg\det G(x)\,dx.
\label{wal}
\neqn
Let us define similarly the lower winding number
$\lw_\al(G)$, by taking $\displaystyle\inf_{y\in\BR}$ and
the corresponding lower limit.
Corollary 1 holds, in particular, for
these characteristics of $G$.
If $n=1$, $G(x)=\exp\big(ik(\sign x)\cdot|x|^\al\big)$,
and $0<\al\le 1$, then
$\uw_\al(G)=\lw{}_\al(G)= k$.


In fact, we could take
instead of $T^{1+\al}$  any function $\rho(T)$ such that
$\rho(T)>0$, $\frac T{\rho(T)}\to0$ as $T\to+\infty$
in the above definitions of generalized winding numbers.

\begin{cor2th3}
Let $G\in \GCarr(\BC_+)$ or
$G\in \GCarr(\BC_-)$,
where $\BC_-=\{z\in\BC:\, \Im z<0\}$. Then for any $\al>0$,
$\uw_\al(G)=\lw{}_\al(G)= 0$.
\end{cor2th3}

Indeed, in both cases $T_G^{-1}=T_{G^{-1}}$, hence
$T_G\in \Phi(H^2_r)$, and we can apply Corollary 1.
\hfill $\Box$

\begin{cor3th3}
Let $G\in \GCbrr(\BR)$, and define
$\uw_1(G)$ by (\ref{wal})
and $\be(G)$ by (\ref{betaG}). Then
$\ds\beta(G)\ge\frac{\uw_1(G)}r$.
\end{cor3th3}

Indeed, if $T_{e^{-i\tau x}G}$ is onto, then it is a
$\Phi_-$ operator, which implies that
\[
\qquad \qquad \qquad \qquad \uw_1(e^{-i\tau x}G)= \uw_1(G) - r\tau \le 0.
\qquad \qquad \qquad \qquad  \Box
\]

We remark that if $G$ is a semi-almost periodic $r\times r$
matrix function
such that $G, G^{-1}\in C^b(\BR)$, then $\det G$ is a scalar semi-almost periodic
function, and $\det G$ has almost periodic representatives
$(\det G)_{\pm \infty}$ at $+\infty$ and $-\infty$, respectively
(see \cite{BKSpi}). These representatives, by the Bohr mean motion theorem,
have the form
$$
(\det G)_{\pm \infty}(x)=e^{\kappa_\pm x}e^{g_\pm(x)},
$$
where $\kappa_\pm $ are mean motions of
$\det G(x)$ at $\pm \infty$ and functions
$g_\pm$ are almost periodic. In this case,
\begin{gather*}
\lw{}_1(G)=\min(\kappa_-, \kappa_+), \qquad
\uw_1(G)=\max(\kappa_-, \kappa_+), \\
\ltw{}_1(G)=\utw{}_1(G)=\frac {\kappa_-+\kappa_+}2.
\end{gather*}
If $r=1$, complete criteria of Fredholmness,
as well as the calculation of the Fredholm index are
known since the work by Sarason \cite{Sar2}.
It follows, in particular, that in this case
$\beta(G)=\max(\kappa_-, \kappa_+)$. So
Corollary 3 of Theorem 3 gives an exact estimate
for the case of scalar almost periodic
functions.

The study of the almost periodic and
semi-almost periodic matrix cases depend on the
existence of some special factorizations of $G$.
If these factorizations exist, then complete
criteria for Fredholmness and formulas for
the index are available, see
\cite{KSpi} and \cite{BKSpi}, Ch. 10 and \S19.6.

\begin{proof}[Proof of Theorem 2]
By Theorem 1, it  only has to be proved that if $f\in\BMOp$,
then
$$
\int_\BR f(x)\eta(x)\,dx\gews -C\|\eta\|_{H^1_\BR}
\qquad
\text{ for all }\eta\in \Pi.
$$
This inequality follows from the Fefferman duality
$H^1_\BR$ -- $\BMOR$ (see \cite{GR})
in the case when $f\in \BMOR$.
Now let $f$ be nondecreasing, and take any function $\eta\in\Pi$.
Suppose that $\supp \eta\subset I$, where
$I$ is a finite interval. Approximate $f$
in $L^\infty(I)$ by
a sequence of nondecreasing step functions
$\{f_n\}$
of the form
$$
f_n=C_n+\sum \al_{nk}\,\chi_{(-\infty,a_{nk}]},
$$
where $C_n, a_{nk}\in\BR$ and $\al_{nk}$ are negative.
 Then
$\int_{\BR}\eta f_n\gews0$ for all $n$,
hence $\int_{\BR}\eta f\gews0$.

We obtain the result by combining these two cases.
\end{proof}

\begin{proof}[Proof of Theorem 3]
Let $\eta_{T,y}(x)=\eta\big(\frac{x-y}T\big)$. Since $\CH(\eta_{T,y})=(\CH\eta)_{T,y}$,
it follows that $\|\eta_{T,y}\|_{H^1_\BR}= T\|\eta\|_{H^1_\BR}$.
So the assertions follow directly from
\linebreak 
Theorem 2.
\end{proof}


\section{Some related questions }
\setcounter{equation}{0}

\nin{\bf Problem 1.} \enspace Give a real variable characterization of classes
$\BMOR^\pm$.

The next two questions are certainly known for specialists
for a long time, however, complete answers are not known.

\nin{\bf Problem 2.} \enspace
1) Let $r=1$, and let
$G\in C(\BR)$, $\arg G\in C_+(\BR)$,  \linebreak 
$\lim _{x\to\pm\infty}\arg G=\pm\infty$.
What additional conditions guarantee that $T_G\in \Phi_+(\Htwor)$?

2) What can be said in this respect for the matrix case $r>1$?

Sufficient conditions for $r=1$
are given in \cite{BG} and in \cite{DybGru},
Theorem 5.10. As it follows from the construction of Lemma 4.9 in
\cite{BG}, there are symbols $G$ of the above type
such that $T_G$ is not semi-Fredholm.
See also \cite{Gov}, Theorem 28.2 and Section 32 for related counter-examples.

The book \cite{DybGru} also contains results about the matrix valued case.
At least for the scalar case, it seems that more complete answers can be found.


 \nin{\bf Problem 3.} \enspace  Suppose that $T_G\in\Phi_+$.
 Can one give some estimates of $\ind T_G$ in terms of some
 explicit real variable characteristics of $\arg\det G$?

\nin{\bf Problem 4.} \enspace
Suppose that $\eta_1, \eta_2\in\Pi$.
When can one assert that $\uw_{\eta_1}(G)\ne+\infty$
implies that $\uw_{\eta_2}(G)\ne+\infty$
for all $G\in\Cbrr(\BR)$ with $|\det G|>\eps>0$ on $\BR$?
Is there a ``universal'' function
$\eta_0\in\Pi$ such that for any  $G$ as above,
$\uw_{\eta_0}(G)\ne+\infty$
implies that $\uw_{\eta}(G)\ne+\infty$
for all $\eta\in\Pi$?

\vskip-.5cm


%


%
%
%

\end{document}